\documentclass[12pt]{amsart}
\usepackage{amssymb,amsmath,amsfonts,amsthm,bbm,mathrsfs,times,graphicx,color,comment,mathabx,enumerate}
\usepackage[utf8]{inputenc}
\usepackage{tikz}
\usetikzlibrary{arrows,shapes}
\usetikzlibrary{fadings}
\usetikzlibrary{intersections}
\usetikzlibrary{decorations.pathreplacing}

\usepackage{scalerel}

\renewcommand{\i}{\mathrm i}

\newcommand{\INF}{{\infty}}

\newcommand{\OM}{\Omega}

\newcommand{\sph}{{{\mathbf S}^ 1}}

\newcommand{\del}{\partial}

\newcommand{\ol}{\overline}
\newcommand{\ds}{\displaystyle}
\newcommand{\dba}{\overline{\partial}}

\newcommand{\BR}{\mathbb{R}}

\newcommand{\bu}{{\bf u}}
\newcommand{\bv}{{\bf v}}

\newcommand{\bg}{{\bf g}}

\newcommand{\bF}{{\bf F}}

\newcommand{\bzero}{\mathbf 0}

\newcommand{\btheta}{\boldsymbol \theta}
\newcommand{\balpha}{\boldsymbol \alpha}
\newcommand{\bbeta}{\boldsymbol \beta}

\newcommand{\B}{\mathcal{B}}

\newtheorem{theorem}{Theorem}[section]
\newtheorem{prop}{Proposition}[section]
\newtheorem{lemma}{Lemma}[section]

\title[Inversion of the exponential X-ray transform of symmetric 2-tensors]
{Inversion of the exponential X-ray transform of symmetric 2-tensors}

\begin{document}
	\date{\today}

	\author{David Omogbhe}
	\address{Faculty of Mathematics, University of Vienna, Oskar-Morgenstern-Platz 1, 1090 Vienna, Austria}
	\email{david.omogbhe@univie.ac.at}

	
	\maketitle
	
	\begin{abstract}
A unique inversion of the exponential X-ray transform of some class of symmetric 2-tensor fields supported in a two-dimensional strictly convex set is presented. The approach to inversion is based on the Cauchy problem for a Beltrami-like equation associated with A-analytic maps.
	\end{abstract}
		
\section{Introduction}
The exponential X-ray transform arises in emission computed tomography. If attenuation of a medium is everywhere constant in a convex set containing the sources of emitted radiation,
then the attenuated X-ray transform reduces to a transform called the exponential X-ray transform. For the exponential X-ray transform of 0-tensors (functions), inversion formulas can be found in \cite{3,4,5}, identification problem in \cite{2,1}, and its range characterization in \cite{9,7,6,8,sadiqScherzerTamasan, D1,sadiqTamasan01,sadiqTamasan02}. Unlike the attenuated X-ray transform of planar 0- and 1-tensors which are fully injective \cite{sg, fn,rg,sadiqTamasan02,d}, the attenuated X-ray transform of symmetric $m$-tensor fields has large kernel, and a vast literature in tensor tomography concerns what part of the symmetric $m$-tensor field is reconstructible from its X-ray transform, see  \cite{r,a,k,a2,i,11,12}.
 
The X-ray transform of symmetric 2-tensors occurs in the linearization of boundary rigidity problem \cite{r}. In this work, we consider a unique inversion of the exponential X-ray transform of some class symmetric 2-tensor  supported in a strictly convex bounded subset in the Euclidean plane. We show that the exponential X-ray transform of symmetric 2-tensors in Theorem\ref{thm3} is invertible in the class of incompressible symmetric 2-tensors  and in Theorem \ref{thm4} in class of trace-free symmetric 2-tensors. The recovery of the symmetric 2-tensor field from its exponential X-ray transform herein is recast as the recovery of a source term in a transport equation from boundary measurements. Our approach to reconstruction is based on the Cauchy problem for Beltrami-like equation associated with $A$-analytic maps in the sense of Bukhgeim \cite{bukhgeimBook}. The  $A$-analytic theory for non-attenuating medium was developed in \cite{bukhgeimBook} and attenuating case was treated in \cite{ABK}.  The Bukhgeim’s approach is based on expanding the solution to the transport equation into its Fourier series with respect to the velocity $\btheta$ and studying the corresponding system of Fourier coefficients by some version of complex analysis. In section 2, we recall some basic properties of $A$-analytic theory and in section 3, we provide the reconstruction methods. \\

For a real symmetric  2-tensor $\bF \in L^1(\mathbb{R}^2;\mathbb{R}^{2\times2}), x\in \mathbb{R}^2:$	

\begin{equation}\label{tensor}
\bF(x) = 
\begin{pmatrix}
f_{11}(x) & f_{12}(x)  \\
f_{12}(x) & f_{22}(x)  
\end{pmatrix},
\end{equation}

and a real number $\mu$ (usually positive), the exponential $X$-ray transform of symmetric 2-tensor field $\bF$ is defined by 

\begin{equation}\label{superpose}
\displaystyle X_\mu\bF(x,\btheta):=\int\limits_{-\infty}^{\infty}\langle\;\bF(x+t\btheta)\btheta,\btheta \rangle\; e^{\mu t}dt,
\end{equation}

where $\langle.,.\rangle$ is the scalar product in $\mathbb{R}^2$, and $\btheta$ is a direction in the unit sphere $\sph.$\\

Throughout, $\bF$  is supported in the interior of a strictly convex  domain $\OM\subset \BR^2$ with $C^{2,\alpha}$ smooth boundary $\partial\Omega,$ for $\alpha > 1/2.$ For any $(x,\btheta)\in \overline{\Omega}\times \sph,$ let $\tau_\pm(x,\btheta)$ denote the distance from $x$ in the $\pm \btheta$ direction to the boundary, and distinguish the endpoints $x_{\btheta}^{\pm} \in \partial \OM$ of the chord $\tau(x,\btheta)= \tau_-(x,\btheta)+\tau_+(x,\btheta) $  in the direction of $\btheta$ passing through $x$ by
\begin{equation*}
x_{\btheta}^\pm = x\pm \tau_\pm(x,\btheta) \btheta.
\end{equation*}	

The incoming (-) respectively outgoing (+) sub-manifolds of the unit bundle restricted to the boundary:

\begin{equation}
\partial \Omega_{\pm} :=\{(\zeta,\btheta)\in \partial\Omega \times\sph:\, \pm\eta(\zeta)\cdot\btheta>0 \}\textcolor{red}{,}
\end{equation}
	with $\eta(\zeta)$ being the outer unit normal at $\zeta \in \partial \Omega.$
\vspace{0.2cm}
	
The exponential $X$-ray transform in \eqref{superpose} is realised as a function on 
$\partial \Omega_+$	 by

\begin{equation}
\displaystyle X_\mu\bF(x,\btheta):=\int\limits_{-\tau(x,\btheta)}^{0} \langle\;\bF(x+t\btheta)\btheta,\btheta \rangle\; e^{\mu t}dt. 
\end{equation}
We approach the reconstruction through its connection with the transport model as follows: The boundary value problem
\begin{subequations}\label{Transportequ}
\begin{eqnarray}
\btheta\cdot\nabla u(x,\btheta) + \mu u(x,\btheta)& = & \langle \bF(x)\btheta,\btheta \;\rangle \quad (x,\btheta)\in \OM\times\sph, \\
u\Big\lvert_{\partial \Omega_-}&=&0,
\end{eqnarray}
\end{subequations}	
has a unique solution in $\Omega\times \sph$, and measurement from the exiting radiation on the boundary $\partial \Omega_+:$	
\begin{equation}\label{data}
u\Big\lvert_{\partial \Omega_+}(x,\btheta)= X_\mu \bF(x,\btheta),\quad (x,\btheta)\in \partial \Omega_+.
\end{equation}
	
	\section{Preliminaries}\label{sec:prelim}

Let $\textbf{w}$ be a planar symmetric $m-$ tensor field, let $k\geq1$ and $m\geq 1$ be integers, $$\delta :H^k(S^m)\longrightarrow H^{k-1}(S^{m-1})$$ is the divergence operator defined by
\begin{align}\label{diveqn}
( \delta \textbf{w} )_{i_1 \cdots i_{m-1}}= \frac{\del w_{i_1  \cdots i_{m-1}j}}{\del x_{j}} = \frac{\del w_{i_1  \cdots i_{m-1}1}}{\del x_{1}}+\frac{\del w_{i_1  \cdots i_{m-1}2}}{\del x_{2}}\textcolor{red}{.}
\end{align}

For symmetric 2-tensor $\bF$ as in \eqref{tensor},
\begin{align}\label{divprop}
\displaystyle\delta (\bF(x)):=\big((\delta\bF)_1,(\delta\bF)_2\big) = \big(\partial_{x_1}f_{11}(x)+\partial_{x_2}f_{12}(x),\partial_{x_1}f_{12}(x) +\partial_{x_2}f_{22}(x)\big).
\end{align}
\vspace{0.2cm}
	
Let $l_1$ and $l_\infty$ be the space of summable and bounded sequences respectively. For $0<\alpha <1$, we consider the  Banach space:
	\begin{equation}\label{spaces}
		\begin{aligned} 
			Y_{\alpha}(\partial \Omega):= \left \{ \bg= ( g_{0}, g_{-1}, g_{-2},...) :\sup_{\xi \in \partial \Omega}\sum_{j=0}^{\INF}  {\langle j \rangle}^2 \lvert g_{-j}(\xi) \rvert < \INF  ~~~\text{and} \;
			\underset{{\substack{{\xi}_1,{\xi}_2 \in \partial \Omega \\
						{\xi}_1 \neq {\xi}_2 } }}{\sup} \sum_{j=0}^{\INF}  \langle j \rangle 
			\frac{\lvert g_{-j}({\xi}_1) - g_{-j}({\xi}_2)\rvert }{|{\xi}_1 - {\xi}_2|^{ \alpha}} < \INF  \right \},
		\end{aligned}
	\end{equation}
where we use for brevity the notation $\langle j \rangle =(1+|j|^2)^{1/2}.$\\

For $z=x_1+\i x_2$, we consider the Cauchy-Riemann operators
	\begin{align}\label{dbar_op} 
		\ol{\del} = \left( \del_{x_{1}}+\i \del_{x_{2}} \right) /2 ,\quad \del = \left( \del_{x_{1}}- \i \del_{x_{2}} \right) /2.
	\end{align}	
	
	A sequence valued map $\OM \ni z\mapsto  \bv(z): = ( v_{0}(z), v_{-1}(z),v_{-2}(z),... )$ in $C(\ol\OM;l_\INF)\cap C^1(\OM;l_\INF)$
	is called {\em $\mathcal{L}^2$-analytic} (in the sense of Bukhgeim), if
	\begin{equation}\label{Aanalytic}
		\ol{\del} \bv (z) + \mathcal{L}^2\del \bv (z) = 0,\quad z\in\OM,
	\end{equation}
	where $\mathcal{L}$ is the left shift operator $\ds \mathcal{L} ( v_{0}, v_{-1}, v_{-2}, \cdots  ) =  ( v_{-1}, v_{-2},  \cdots ), $ and
	$\mathcal{L}^2=\mathcal{L}\circ \mathcal{L}$.

	Bukhgeim's original  theory \cite{bukhgeimBook}  shows that solutions of \eqref{Aanalytic},  satisfy a Cauchy-like integral formula,
	\begin{align}\label{Analytic}
		\bv (z) = \B [\bv \lvert_{\partial \Omega}](z), \quad  z\in\OM,
	\end{align} where $\B$ is 
	the Bukhgeim-Cauchy operator  acting on $\bv \lvert_{\partial \Omega}$. We use the formula in \cite{finch}, where
	$\B$ is defined component-wise for $n\geq 0$ by
	\begin{equation} \label{BukhgeimCauchyFormula}
		\begin{aligned} 
			(\B \bv)_{-n}(z) &:= \frac{1}{2\pi \i} \int\limits_{\partial \Omega}
			\frac{ v_{-n}(\zeta)}{\zeta-z}d\zeta  
			+ \frac{1}{2\pi \i}\int\limits_{\partial \Omega} \left \{ \frac{d\zeta}{\zeta-z}-\frac{d \ol{\zeta}}{\ol{\zeta}-\ol{z}} \right \} \sum_{j=1}^{\infty}  
			v_{-n-2j}(\zeta)
			\left( \frac{\ol{\zeta}-\ol{z}}{\zeta-z} \right) ^{j},\; z\in\OM.
		\end{aligned}
	\end{equation}
	

	\begin{theorem}\cite[Proposition 2.3]{sadiqTamasan02}.\label{BukhgeimCauchyThm}
		Let $\B$ be the Bukhgeim-Cauchy operator in \eqref{BukhgeimCauchyFormula}. 
		
		If $\bg = ( g_{0}, g_{-1}, g_{-2},...) \in Y_{\alpha}(\partial \Omega)$ for $\alpha>1/2$, then $ \bv :=\B \bg\in C^{1,\alpha}(\OM;l_1)\cap C^{\alpha}(\ol \OM;l_1)\cap C^2(\OM;l_\infty)$ , and $\bv$ is $\mathcal{L}^2$-analytic in $\OM$. 
	\end{theorem}

In addition to $\mathcal{L}^2$-analytic maps,  consists of the one-to-one relation between solutions
	$ \bu: = ( u_{0}, u_{-1},u_{-2},... ) $
	satisfying
	\begin{align}\label{beltrami}
		\dba\bu +\mathcal{L}^2 \del\bu+ a\mathcal{L}\bu = \bzero,
	\end{align}
	 
	and the $\mathcal{L}^2$-analytic map $\bv = ( v_{0}, v_{-1},v_{-2},... )$ satisfying \eqref{Aanalytic} via a special function $h$, see \cite[Lemma 4.2]{sadiqTamasan02} for details.
	The function $h$ is defined as 
	\begin{align}\label{hDefn}
		h(z,\btheta) := \ds \int\limits_{0}^{\INF} a(z+t\btheta)dt -\frac{1}{2} \left( I - \i H \right) Ra(z\cdot \btheta^{\perp}, \btheta^{\perp}),
	\end{align} 
	
where $\btheta^\perp$ is  the counter-clockwise rotation of $\btheta$ by $\pi/2$,
	$Ra(s, \btheta^{\perp}) = \ds \int\limits_{-\INF}^{\INF} a\left( s \btheta^{\perp} +t \btheta \right)dt$ is the Radon transform in $\BR^2$  of  $a$, and $H$ is the classical Hilbert transform where
	
\begin{align*}
HRa(z\cdot \btheta^{\perp}, \btheta^{\perp}) = H[Ra(\cdot, \btheta^{\perp})](z\cdot \btheta^{\perp}) = \ds \frac{1}{\pi}\int\limits_{-\INF}^{\INF} \frac{Ra(t, \btheta^{\perp})}{(z\cdot \btheta^{\perp})-t}dt.
\end{align*}	
	 
The function $h$ appeared first in \cite{nattererBook} and enjoys the crucial property of having vanishing negative Fourier modes yielding the expansions
	
	\begin{align}\label{ehEq}
		e^{- h(z,\btheta)} := \sum_{k=0}^{\INF} \alpha_{k}(z) e^{\i k\varphi}, \quad e^{h(z,\btheta)} := \sum_{k=0}^{\INF} \beta_{k}(z) e^{\i k\varphi}, \quad (z, \btheta) \in \ol\OM \times \sph.
	\end{align}
	
	Using the Fourier coefficients of  $e^{\pm h}$,  define the operators $e^{\pm G} \bu $ component-wise for each $n \leq 0$, by 
	\begin{align}\label{eGop}
		(e^{-G} \bu )_n &= (\balpha \ast \bu)_n = \sum_{k=0}^{\infty}\alpha_{k} u_{n-k}, \quad \text{and} \quad 
		(e^{G} \bu )_n = (\bbeta \ast \bu)_n = \sum_{k=0}^{\infty}\beta_{k} u_{n-k}, \quad \text{where} \\ \nonumber
		&\ol \OM \ni z\mapsto \balpha(z) := ( \alpha_{0}(z), \alpha_{1}(z),  ... , ),  
		\quad \ol \OM \ni z\mapsto \bbeta(z) := ( \beta_{0}(z), \beta_{1}(z),  ... , ). 
	\end{align}


We remark that in this work, the function $a$ will be taken as the constant attenuation $\mu$ defined inside the convex set containing the support of the symmetric 2-tensor field $\bF.$	
	
	\begin{prop}\cite[Proposition 5.2]{sadiqTamasan01}\label{eGprop}
Let $e^{\pm G}$ be operators as defined in \eqref{eGop}, and Let  $\alpha >1/2$. If $\bg \in Y_{\alpha}(\partial \Omega)$ then $e^{\pm G}\bg \in Y_{\alpha}(\partial \Omega) $

	\end{prop}
	\begin{lemma}\cite[Lemma 4.2]{sadiqTamasan02}\label{beltrami_reduction}
		Let $e^{\pm G}$ be operators as defined in \eqref{eGop}. 
		
		(i) If $\bu \in C^1(\OM, l_1)$ solves $\ds \dba \bu + \mathcal{L}^2 \del \bu+ \mu \mathcal{L}\bu = \bzero$, then  $\ds \bv= e^{-G} \bu \in C^1(\OM, l_1)$ solves $\dba \bv + \mathcal{L}^2\del \bv =\bzero$.
		
		(ii) Conversely, if $\bv \in C^1(\OM, l_1)$ solves $\dba \bv + \mathcal{L}^2\del \bv =\bzero$, then $\ds \bu= e^{G} \bv \in C^1(\OM, l_1)$ solves $\ds \dba \bu +\mathcal{L}^2 \del \bu+ \mu \mathcal{L}\bu = \bzero$.
	\end{lemma}
\vspace{0.2cm}

	\section{inversion of the exponential X-ray transform of symmetric 2-tensors}\label{sec:reconstruct}
For a real valued symmetric 2-tensor field $\bF$ in \eqref{tensor}, and $\btheta = (\cos\varphi, \sin \varphi),$ we have	
	\begin{equation}\label{Fourierexp}
	\displaystyle\langle ~\bF(x)\btheta,\btheta\rangle = f_0(z) +\overline{f_2(z)}e^{2i\varphi} + f_2(z)e^{-2i\varphi}\quad 
	\end{equation}
	where	
\begin{equation}\label{Fouriercoeff}
	f_0(z) =\frac{f_{11}(x)+f_{22}(x)}{2},\quad f_2(z) =\frac{f_{11}(x)-f_{22}(x)}{4}+i\frac{f_{12}(x)}{2}\ .
	\end{equation}

From \eqref{Fourierexp}, the transport equation \ref{Transportequ}(a) becomes
	\begin{equation}\label{trans2}
	\btheta \cdot \nabla u(z,\btheta)+ \mu u(z,\btheta) = f_0(z) +\overline{f_2(z)}e^{2i\varphi} + f_2(z)e^{-2i\varphi}.
	\end{equation}
\vspace{0.2cm}
For $z = x_1+ix_2 \in \Omega$, let $u(z,\btheta) = \sum_{-\infty}^{\infty} u_{n}(z) e^{in\varphi}$ be the formal Fourier series representation of the solution of \eqref{trans2} in the angular variable $\btheta=(\cos\varphi,\sin\varphi)$. Since $u$ is real valued, the Fourier coefficients $\{u_{n}\}$  occurs in complex-conjugate pairs $u_{-n}=\ol{u_n}$. For	the derivatives  $\del,\dba$  in the spatial variable as in \eqref{dbar_op}, the advection operator $\btheta \cdot\nabla$ in \eqref{trans2}  becomes  $\btheta \cdot\nabla=e^{-\i \varphi}\dba + e^{\i \varphi}\del$.  By identifying the Fourier coefficients of the same order, the non-positive Fourier coefficients\\
	
\begin{align}\label{boldu}
		\OM \ni z\mapsto  \bu(z)&: = ( u_{0}(z), u_{-1}(z),u_{-2}(z),... )
	\end{align}
of $u$ satisfy:	
		
 
	\begin{align}\label{l1}
		\overline{\del} u_{1}(z)+\del u_{-1}(z)+ \mu u_0(z) &=f_0(z), \\ \label{l2}
		\overline{\del} u_{0}(z)+\del u_{-2}(z)+  \mu u_{-1}(z) &= 0, \\ \label{l3}
		\overline{\del} u_{-1}(z)+\del u_{-3}(z)+ \mu u_{-2}(z) &= f_2(z), \\
		\dba u_{-n}(z) +\del u_{-n-2}(z)+ \mu u_{-n-1}(z)  &=0,\label{l4}\qquad n\geq 2,
	\end{align}
where $f_0, f_2$ are define in \eqref{Fouriercoeff}. 

\begin{lemma}\label{k1}
Let $\OM\subset\mathbb{R}^2$ be a strictly convex bounded domain and $\partial \Omega$ be its boundary. Consider the boundary value problem \eqref{Transportequ} for some known positive real number $\mu$.
	If the unknown real valued symmetric 2-tensor $\bF$ is in $W^{3,p}( \OM; \BR^{2\times2})$ with $p > 4$, 
	then  the data $u\lvert_{\partial \Omega_+}$ uniquely determine in $\OM$, the Fourier coefficients $$\mathcal{L}^2 \bu = ( u_{-2},u_{-3},u_{-4},...) \in C^{1,\alpha}(\OM;l_1)\cap C^{\alpha}(\ol \OM;l_1)\cap C^2(\OM;l_\infty)$$ 
of the solution $u$ to the boundary value problem \eqref{Transportequ}.	
\end{lemma}

\begin{proof}
Let $u$ be the solution of the boundary value problem \eqref{Transportequ} in $\OM$, and let $\bu= ( u_{0}, u_{-1}, u_{-2}, ... ) $ be the sequence of its non-positive Fourier coefficients.
Let $\bg= ( g_{0}, g_{-1}, g_{-2}, ... ) $  be the sequence valued map of its non-positive Fourier coefficients of the data

$$
g(z,\btheta):=u\Big\lvert_{\partial \OM \times\sph}=
\begin{cases}
X_\mu\bF(z,\btheta),\quad (z,\btheta)\in {\partial \OM}_+ \\
0,~~\qquad \qquad (z,\btheta)\in {\partial \OM}_- ,
\end{cases}
$$

where

$$g_{-n}(z) =\frac{1}{\pi}\int\limits_0^{2\pi} g(z,\btheta)e^{in\varphi}d\varphi, \quad n\geq 0, \quad g_{n}(z) = \overline{g_{-n}(z)}.$$

Since the  symmetric 2-tensor $\bF$ is $W^{3,p}( \OM; \BR^{2\times2})$-regular with $p > 4$, then the anisotropic source $\displaystyle f(x,\btheta):= \langle \bF(x)\btheta, \btheta\rangle$ belong to $W^{3,p}(\OM \times \sph)$ with $p > 4$. By applying \cite[Theorem A.1]{D2}, we have $u \in W^{3,p}(\Omega\times\sph), \,p > 4$. Moreover, by the Sobolev embedding \cite{adam},
	$$ W^{3,p}(\Omega\times\sph) \subset C^{2,\alpha}(\ol \OM\times\sph) \quad \mbox{with} \quad 
	\alpha =1-\frac{2}{p}>\frac{1}{2},$$ thus  $u\in C^{2,\alpha}(\ol \OM\times\sph)$,
	and by \cite[Proposition 4.1 (ii)]{sadiqTamasan01}, the sequence valued map $\bg\in Y_{\alpha}(\partial\Omega)$. 
	We note from \eqref{l4} that the shifted sequence valued map $\mathcal{L}^2\bu=( u_{-2},u_{-3},u_{-4},...)$ solves 
	\begin{align}\label{Lmu_beltrami}
		\dba \mathcal{L}^2 \bu(z) +\mathcal{L}^2 \del \mathcal{L}^2 \bu(z)+ \mu \mathcal{L}(\mathcal{L}^2\bu(z)) = \bzero,\quad z\in \OM.
	\end{align}
Let $\bv := e^{-G} \mathcal{L}^2\bu,$ by Lemma \ref{beltrami_reduction}, the sequence $\bv$ is $\mathcal{L}^2$ analytic.	
From the boundary data $\bg \in Y_{\alpha}(\partial\Omega)$ and Proposition \ref{eGprop}, we determines on the boundary $${\bv}\vert_{\partial \OM} := e^{-G} \mathcal{L}^2 \bg = \mathcal{L}^2 e^{-G}  \bg \in Y_{\alpha}(\partial\Omega),$$ where the operators $e^{\pm G}$ commute with the left translation $\mathcal{L}$.
	
From ${\bv}\vert_{\partial \OM} \in Y_{\alpha}(\partial\Omega) $, we use the Bukhgeim-Cauchy Integral formula \eqref{BukhgeimCauchyFormula} to 
	construct the sequence valued map $\bv $ inside $\OM$
	$$\bv(z) = \B[{\bv}\rvert_{\partial\OM}](z), \quad z\in \OM,$$
satisfying $\dba \bv + \mathcal{L}^2\del \bv =\bzero.$ 
By Theorem \ref{BukhgeimCauchyThm}
	the constructed $\mathcal{L}^2$-analytic sequence valued  map $\bv$ in $\OM$ is in $C^{1,\alpha}(\OM;l_1)\cap C^{\alpha}(\ol \OM;l_1)\cap C^2(\OM;l_\infty).$ We use convolution in Lemma \ref{beltrami_reduction} to determine sequence $\ds \mathcal{L}^2 \bu := e^{G}\bv $, i.e
\begin{equation}\label{conseq1}
( u_{-2},u_{-3},u_{-4},...)\in C^{1,\alpha}(\OM;l_1)\cap C^{\alpha}(\ol \OM;l_1)\cap C^2(\OM;l_\infty)
\end{equation}	
inside $\OM$ satisfying \eqref{l4}. \qed
\vspace{0.4cm}
  
The theorem below shows that exponential X-ray transform of symmetric 2-tensors is invertible for incompressible symmetric 2-tensors.

\begin{theorem}\label{thm3}
	Let $\OM\subset\mathbb{R}^2$ be a strictly convex bounded domain and $\partial \Omega$ be its boundary. Consider the boundary value problem \eqref{Transportequ} for some known positive real number $\mu$. 
	If the unknown real valued incompressible symmetric 2-tensor $\bF$ is $W^{3,p}( \OM; \BR^{2\times2})$ regular with $p > 4$, 
	then  the data $u\lvert_{\partial \Omega_+}$ uniquely determines the tensor $\bF$ in $\OM.$
\end{theorem}

\begin{proof}
Since $\bF$ is $W^{3,p}( \OM; \BR^{2\times2})$ regular with $p > 4$, and $\mu$ positive real number, then from Lemma\ref{k1}, the boundary measurement $u\lvert_{\partial \Omega_+}$ uniquely determines the Fourier coefficients $( u_{-2},u_{-3},u_{-4},...)\in C^{1,\alpha}(\OM;l_1)\cap C^{\alpha}(\ol \OM;l_1)\cap C^2(\OM;l_\infty)$ ~~inside $\OM$ satisfying \eqref{l4}.\\

By multiplying both sides of \eqref{l3} by $4\partial$,  we have
	  
	 \begin{align}\label{BVP3}
\Delta u_{-1}(z) +4\partial^2 u_{-3}(z)+ 4 \mu \partial u_{-2}(z) &= 4\partial f_2(z)\nonumber\\ \nonumber
&=\big(\partial_{x_1}f_{11}(x)+\partial_{x_2}f_{12}(x)\big)-\frac{1}{2}\partial_{x_1}\left(f_{11}(x)+f_{22}(x)\right)\\ \nonumber&+ i\big(\partial_{x_1}f_{12}(x)+\partial_{x_2}f_{22}(x)\big)-\frac{1}{2}i\partial_{x_2}\left(f_{11}(x)+f_{22}(x)\right)\\&
=(\delta\bF)_1+i(\delta\bF)_2-2\overline{\partial}f_0(z).
	 \end{align}
	 
Since $\bF$ is incompressible i.e $\delta \bF = \bzero$ in \eqref{divprop}, then equation \eqref{BVP3} becomes

\begin{align}\label{k2}
\Delta u_{-1}(z) +4\partial^2 u_{-3}(z)+ 4 \mu \partial u_{-2}(z) = -2\overline{\partial}f_0(z)
\end{align} 

Substituting \eqref{l1} into equation \eqref{k2} gives 

\begin{align}\label{p21}
\Delta u_{-1} +4\partial^2 u_{-3}+4 \mu \partial(u_{-2}) = -4~\overline{\partial}(\mathbb{R}e~\partial u_{-1})-2 \mu \overline{\partial} u_0.
\end{align}

From \eqref{l2}, substituting for $\overline{\partial} u_0$ into \eqref{p21} yields,

\begin{align}\label{p2}
\Delta u_{-1}+ 4~\overline{\partial}(\mathbb{R}e~\partial u_{-1})- 2\mu^2 u_{-1}  = -4\partial^2 u_{-3} -2\mu {\partial} u_{-2}.
\end{align}

Equation \eqref{p2} together with the boundary value yields the boundary value problem
	
	\begin{subequations} \label{p3}
		\begin{eqnarray} \label{Poisson_UM1}
	\Delta u_{-1}(z)+ 4~\overline{\partial}(\mathbb{R}e~\partial u_{-1}(z))- 2\mu^2 u_{-1}(z) & = & -4\partial^2 u_{-3}(z) -2\mu {\partial} u_{-2}(z) ,\quad z\in\Omega,\\  
			\label{UM_Gam1} 
u_{-1}(z)\Big\lvert_{\partial \Omega} &=& g_{-1}(z),
		\end{eqnarray} 
	\end{subequations}
where the right hand side of \eqref{p3} is known. By Theorem \ref{P1}, the BVP \eqref{p3} is uniquely solvable. Thus $u_{-1}$ is uniquely recovered in $\Omega$ by solving the Dirichlet BVP \eqref{p3}.\\

Equation \eqref{l2} together with the boundary value yields the boundary value problem

\begin{subequations} \label{p4}
		\begin{eqnarray} \label{Poisson_UM1}
	\ol{\del} u_{0}(z) & = & -\partial u_{-2}(z)- \mu u_{-1}(z),\quad z\in\Omega,\\  
			\label{UM_Gam1} 
u_{0}(z)\Big\lvert_{\partial \Omega} &=& g_{0}(z),
		\end{eqnarray} 
	\end{subequations}
	
where the right hand side of \eqref{p4} is known. We uniquely solve \eqref{p4} as a Cauchy problem for the $\ol \del$ -equation with boundary value on $\partial \OM$

\begin{subequations} \label{pk}
		\begin{eqnarray} \label{Poisson_UM1}
	\ol{\del} u_{0} & = & \Psi ,\quad \mbox{in} \quad \Omega,\\  
			\label{UM_Gam1} 
u_{0} &=& g_{0}, \quad \mbox{on} \quad \partial \Omega,
		\end{eqnarray} 
	\end{subequations}
via the Cauchy-Pompieu formula \cite{Vekua}:

\begin{align}\label{u0}
u_0(z) = \frac{1}{2\pi i}\int\limits_{\partial \Omega}\frac{g_0(\zeta)}{\zeta-z}d\zeta -\frac{1}{\pi} \iint\limits_{\Omega} \frac{\Psi(\zeta)}{\zeta - z} d\xi d\eta,\quad \zeta = \xi +i\eta,\quad \Psi =  -\partial u_{-2}- \mu u_{-1}, \quad z\in \OM.
\end{align}


From \eqref{Fouriercoeff}, we determine in $\Omega$ the symmetric 2-tensor:

\begin{equation*}\label{a}
\bF(x) = 
\begin{pmatrix}
f_0 + 2~\mathbb{R}e f_2 & 2~\mathbb{I}m f_2  \\
2~\mathbb{I}m f_2 & f_0 - 2~\mathbb{R}e f_2  
\end{pmatrix},
\end{equation*}
where from \eqref{l1}, and \eqref{l3},
\begin{equation*}\label{b}
f_0 = 2~\mathbb{R}e\partial (u_{-1}(z))+\mu u_0(z) \quad \mbox{and} \quad f_2 =\overline{\partial}(u_{-1}(z))+\partial u_{-3}(z)+ \mu u_{-2}(z),
\end{equation*}
and the Fourier coefficients are obtained from equation  \eqref{conseq1}, \eqref{p3}  and \eqref{u0}.
\end{proof}

Next, we show that the exponential X-ray transform of symmetric 2-tensors is invertible for trace-free symmetric 2-tensors.

\vspace{0.5cm}

\begin{theorem}\label{thm4}
	Let $\OM\subset\mathbb{R}^2$ be a strictly convex bounded domain, and $\partial \Omega$ be its boundary. Consider the boundary value problem \eqref{Transportequ} for some known positive real number $\mu$. If the unknown real valued  trace-free symmetric 2-tensor $\bF$ is $ W^{3,p}( \OM; \BR^{2\times2})$ regular with $p > 4$, 
then the $u\lvert_{\partial \Omega_+}$ uniquely determines the tensor $\bF$ in $\Omega.$ 
\end{theorem}

\begin{proof}
Since $\bF$ is $W^{3,p}( \OM; \BR^{2\times2})$ regular with $p > 4$, and $\mu$ positive real number, then from Lemma \ref{k1}, the boundary measurement $u\lvert_{\partial \Omega_+}$ uniquely determines the Fourier coefficients $( u_{-2},u_{-3},u_{-4},...)\in C^{1,\alpha}(\OM;l_1)\cap C^{\alpha}(\ol \OM;l_1)\cap C^2(\OM;l_\infty)$ ~~inside $\OM$ satisfying \eqref{l4}.\\

The 2-tensor $\bF$ being is trace-free implies $f_{11}+f_{22}=0$,  then from \eqref{Fouriercoeff}, we have

\begin{align*}
f_0= \frac{f_{11}+f_{22}}{2} = 0,
\end{align*}	

$$f_2 = \frac{(f_{11}+f_{22})-2f_{22}}{4}+i\frac{f_{12}}{2} = \frac{2f_{11}-(f_{11}+f_{22})}{4}+i\frac{f_{12}}{2},$$ 
\qquad thus
\begin{align}\label{f2mode}
f_0=0, \quad \mbox{and} \quad f_2 = -\frac{f_{22}}{2} + i\frac{f_{12}}{2} = \frac{f_{11}}{2} + i\frac{f_{12}}{2}.
\end{align}
\phantom\qedhere
\end{proof}
	
By applying $4 \del $ to \eqref{l2}, yields

\begin{align}\label{u_0coeff}
\Delta u_0 +4\mu(\partial u_{-1}) = -4 \partial^2 u_{-2},
\end{align}

where $u_0$ is real valued. From \eqref{l1}, and $f_0=0$ by $\bF$ being trace-free, we have the equation

\begin{align}\label{v2}
2~\mathbb{R}e(\partial u_{-1})= -\mu u_0.
\end{align}

Substituting \eqref{v2} into the real part of \eqref{u_0coeff}, we have
\begin{align}\label{t1}
\Delta u_0 -2 \mu^2 u_0 = -4~\BR e~\partial^2 u_{-2}.
\end{align} 

Equation \eqref{t1} together with the boundary value, yields the boundary value problem
	
	\begin{subequations} \label{BVP5}
		\begin{eqnarray} \label{Poisson_UM}
			\Delta u_0(z) -2 \mu^2 u_0(z) &=& -4~\BR e~\partial^2 u_{-2}(z),\quad z \in \Omega,\\  
			\label{UM_Gam} 
			u_0(z) \big\lvert_{\partial \Omega} &=& g_{0}(z),
		\end{eqnarray} 
	\end{subequations}
where the right hand side of \eqref{BVP5} is known.	Thus $u_0$ is uniquely recovered  in $\Omega$ by solving the Dirichlet problem for the screened Poisson equation \eqref{BVP5}.\\
	
From the Imaginary part of 	\eqref{u_0coeff}, we have

\begin{align}\label{v3}
\mathbb{I}m~(\partial u_{-1}) = -\frac{1}{\mu}~\mathbb{I}m~(\partial^2 u_{-2}).
\end{align}

From \eqref{v2} and \eqref{v3} we have

\begin{align}\label{t2}
\partial u_{-1} = -\frac{\mu}{2}u_0 -\frac{i}{\mu}~\mathbb{I}m~(\partial^2 u_{-2}).
\end{align}

Equation \eqref{t2} together the boundary value, yields the boundary value problem	

\begin{subequations} \label{v4}
		\begin{eqnarray} 
		\partial u_{-1}(z) &=& -\frac{\mu}{2}u_0(z) -\frac{i}{\mu}~\mathbb{I}m~(\partial^2 u_{-2}(z)),\quad z \in \Omega\\  
			u_{-1}(z) \big\lvert_{\partial \Omega} &=& g_{-1}(z),
		\end{eqnarray} 
	\end{subequations}
	
where the right hand side of \eqref{v4} is known. We uniquely solve \eqref{v4} as a Cauchy problem for the $\del$ -equation with boundary data on $\partial \OM$

\begin{subequations} \label{pk}
		\begin{eqnarray} \label{Poisson_UM1}
	{\del} u_{-1} & = & \Psi ,\quad \mbox{in} \quad \Omega,\\  
			\label{UM_Gam1} 
u_{-1} &=& g_{-1}, \quad \mbox{on} \quad \partial \Omega \textcolor{red}{,}
		\end{eqnarray} 
	\end{subequations}
via the Cauchy-Pompieu formula \cite{Vekua}:

\begin{align}\label{t4}
u_{-1}(z) = -\frac{1}{2\pi i}\int\limits_{\partial \Omega}\frac{g_{-1}(\zeta)}{\ol\zeta- \ol z}d\ol\zeta -\frac{1}{\pi} \iint\limits_{\Omega} \frac{\Psi(\zeta)}{\ol\zeta -\ol z} d\xi d\eta,\quad \zeta = \xi +i\eta,\quad \Psi = -\frac{\mu}{2}u_0 -\frac{i}{\mu}~\mathbb{I}m~(\partial^2 u_{-2}), \quad z\in \OM.
\end{align}	

From \eqref{f2mode}, the symmetric 2-tensor is given in $\Omega$ by	
\begin{equation}
\bF(x) = 
\begin{pmatrix}
2\mathbb{R}e f_2 & 2~\mathbb{I}m f_2  \\
2~\mathbb{I}m f_2 & -2\mathbb{R}e f_2  
\end{pmatrix}, \quad \mbox{where}~~\mbox{from}~~ \eqref{l3},~~f_2 =\overline{\partial}u_{-1} +\partial u_{-3} + \mu u_{-2},
\end{equation}

and the Fourier coefficients are obtained from
equations \eqref{conseq1}, \eqref{BVP5}, and \eqref{t4}.
\end{proof}


\appendix
\section{}
 The theorem below follows from the generalization of \cite[Theorem A.2]{D2}. Throughout, for repeated indices, we mean Einstein summation convention for sum from 1 to 2.

\begin{theorem}\label{P1}
Let $\OM\subset \BR^2$ strictly convex domain with $C^{2,\alpha}$ smooth boundary $\partial\Omega.$ For a continuous right hand side, the boundary value problem \eqref{BVP4}
is uniquely solvable.

\begin{subequations}\label{BVP4}
\begin{eqnarray}
\Delta u(z) + C_1~\overline{\partial}\big(\mathbb{R}e~\partial u(z)\big)+ C_2 u (z) & =& f(z), \quad z\in \Omega, \\
u(z)\Big\lvert_{\partial\Omega}& =& g(z),
\end{eqnarray}
\end{subequations}
\end{theorem}

where the constant $C_1\geq 0$, and the constant  $C_2\leq 0.$

\begin{proof}
Let $u = u^1 +i u^2,$ then the real and imaginary parts of \eqref{BVP4} satisfies the BVP

\begin{align}\label{system}
P(D)\bar{u} = \tilde{f}, \quad \bar{u}|_{\partial \Omega} = \tilde{g}
\end{align}
where

\begin{equation*}
P(D)= 
\begin{pmatrix}
\Delta + b~ {\partial}^2_{x_1}+C_2 & b~ {\partial}_{x_1}{\partial}_{x_2}  \\
 b~ {\partial}_{x_1}{\partial}_{x_2} &   \Delta + b~ {\partial}^2_{x_2}+C_2
\end{pmatrix},\quad b = C_1/4,
\end{equation*}
\begin{align*}
 \quad \bar{u} = (u^1,u^2)^\top, \quad \tilde{f}=(\mathbb{R}e~f(z), \mathbb{I}m~f(z))^\top \quad \mbox{and}~~ \tilde{g}= (\mathbb{R}e~g(z),\mathbb{I}m~g(z))^\top.
\end{align*}

To show existence and uniqueness of solution, it suffices to show (see \cite{r, taylor,a3}) that $P(D)$ is uniformly elliptic with zero kernel and zero co-kernel.\\\\
We first prove the uniform ellipticity of $P(D)$. For $\zeta =(\zeta_1,\zeta_2) \in {\mathbb{R}}^2 \backslash \{\bzero\} $, the symbol of $P(D)$ is given by

\begin{equation*}
P(\zeta)= 
\begin{pmatrix}
|\zeta|^2 + b ~{\zeta}^2_1+C_2 & b~{\zeta}_1{\zeta}_2  \\
b~{\zeta}_1{\zeta}_2 & |\zeta|^2 + b~{\zeta}^2_2+C_2
\end{pmatrix},
\end{equation*}

and the principal symbol of $P(D)$ is given by

\begin{equation*}
P_L(\zeta)= 
\begin{pmatrix}
|\zeta|^2 + b ~{\zeta}^2_1 & b~{\zeta}_1{\zeta}_2  \\
b~{\zeta}_1{\zeta}_2 & |\zeta|^2 + b~{\zeta}^2_2
\end{pmatrix}.
\end{equation*}

Let $\rho = (\rho_1,\rho_2)\in {\mathbb{R}}^2 \backslash \{\bzero\}$, then $$\rho P_L(\zeta) \rho^T = |\rho|^2|\zeta|^2 + b (\rho_1\zeta_1+\rho_2\zeta_2)^2 \geq|\rho|^2|\zeta|^2 .$$ 

This proves that $P_L(\zeta)$ is positive definite and thus $P(D)$ is uniformly elliptic.\\

Next, for zero right hand side in \eqref{system}, we need  to show that zero vector is the only solution .\\\\
$P(D)\bar{u}= \bzero,$  implies

\begin{subequations}\label{systemBVP}
\begin{eqnarray}
\Delta u^1 + b~{\partial}_{x_1}\nabla\cdot\bar{u}+ C_2~ u^1 &=&0\label{systemBVPa},\\
\Delta u^2 + b~{\partial}_{x_2}\nabla\cdot\bar{u} + C_2~ u^2& =&0\label{systemBVPb}.
\end{eqnarray}
\end{subequations}

Multiplying equation \eqref{systemBVPa} by $u^1$ and \eqref{systemBVPb} by $u^2$ and integrating over $\Omega$ where $\bar{u}|_{\partial \Omega} = \bzero,$ we have

\begin{subequations}\label{integral}
\begin{eqnarray}
-\int\limits_{\Omega}|\nabla u^1|^2 dx - b~\int\limits_{\Omega}\nabla \cdot \bar{u}\partial_{x_1}u^1 dx + C_2~\int\limits_{\Omega}(u^1)^2 dx &=&0\label{integrala},\\
-\int\limits_{\Omega}|\nabla u^2|^2 dx - b~\int\limits_{\Omega}\nabla \cdot \bar{u}\partial_{x_2}u^2 dx + C_2~\int\limits_{\Omega}(u^2)^2 dx& =&0.\label{integralb}
\end{eqnarray}
\end{subequations}
Adding equation \eqref{integrala} and \eqref{integralb}, we have

\begin{align}\label{zp}
\int\limits_{\Omega}\Big(|\nabla u^1|^2+|\nabla u^2|^2 + b~(\nabla \cdot \bar{u})^2- C_2|\bar{u}|^2 \Big)dx = 0.
\end{align}

Since $b = C_1/4\geq 0$  and  $C_2\leq 0,$ it is easy to see that $\bar{u} = \bzero$  is the only solution to \eqref{zp}, hence Ker$(P(D)) = \{\bzero\}.$\\

Finally, we show that the co-kernel of $P(D)$ is also zero.\\ Let $\bar{\omega} = (\omega^1,\omega^2)\in\big(\mbox{Image}(P(D)\big)^\perp$. i.e we consider $\bar{\omega}$ such that

\begin{align}\label{product}
\langle\bar{\omega}, P(D)\bar{u}\rangle_{L^2(\Omega)}=0\quad \mbox{for all}\quad \bar{u}\in C^{\infty}(\Omega)\quad \mbox{with}\quad \bar{u}|_{\partial \Omega} = \bzero.
\end{align}

Since $P(D)$ is self adjoint, it implies

\begin{align}\label{product2}
\langle P(D)\bar{\omega}, \bar{u}\rangle_{L^2(\Omega)}=0\quad \mbox{for all}\quad \bar{u}\in C^{\infty}(\Omega)\quad \mbox{with}\quad \bar{u}|_{\partial \Omega} = \bzero.
\end{align}

Since $\bar{u}\in C^{\infty}_0(\Omega)$ is arbitrary, we have

\begin{align}\label{cokernel}
P(D)\bar{\omega} = \bzero.
\end{align}

From \eqref{product}, we have
\begin{align}
\langle\bar{\omega}, P(D)\bar{u}\rangle_{L^2(\Omega)}=\int\limits_{\Omega
}\big(\omega^i\Delta u^i + b~\omega^i\partial_{x_i}\nabla\cdot\bar{u} + C_2~\omega^i u^i \big)dx = 0.
\end{align}
From Green's second formula, where $\bar{u}\in C^{\infty}_0(\Omega)$ we have

\begin{align}\label{integral1}
\int\limits_{\Omega}\omega^i\Delta u^i dx =\int\limits_{\Omega}u^i\Delta \omega^i dx +\int\limits_{\partial\Omega}\omega^i\nabla u^i\cdot\eta ds,
\end{align}
and integration by parts yields

\begin{align}\label{integral2}
\int\limits_{\Omega}\omega^i\partial_{x_i}\nabla\cdot\bar{u}dx = \int\limits_{\Omega}u^i\partial_{x_i}\nabla\cdot\bar{\omega}dx + \int\limits_{\partial\Omega}\bar{\omega}\cdot \eta \nabla\cdot \bar{u} ds.
\end{align}
From equation \eqref{integral2}, we have
\begin{align}\label{L1}
b\int\limits_{\Omega}\omega^i\partial_{x_i}\nabla\cdot\bar{u}dx + \int\limits_{\Omega} C_2~\omega^i u^i dx = b\int\limits_{\Omega}u^i\partial_{x_i}\nabla\cdot\bar{\omega}dx + \int\limits_{\Omega} C_2~\omega^i u^i dx + b\int\limits_{\partial\Omega}\bar{\omega}\cdot \eta \nabla\cdot \bar{u} ds .
\end{align}

Adding \eqref{integral1} and \eqref{L1} we have
 
 \begin{align}\label{kp}
\langle\bar{\omega}, P(D)\bar{u}\rangle_{L^2(\Omega)} = \langle P(D)\bar{\omega}, \bar{u}\rangle_{L^2(\Omega)} + \langle\bar{\omega}, \mathbb{G}(\partial_\eta)\bar{u}\rangle_{L^2(\partial\Omega)}
 \end{align}
where $\mathbb{G}(\partial_\eta)$ is the boundary operator defined by

\begin{align}
\mathbb{G}(\partial_\eta)\bar{u} = (\nabla u^1\cdot\eta,\nabla u^2\cdot\eta) + b~(\nabla\cdot\bar{u})\eta.
\end{align}

From \eqref{product} and \eqref{product2}, equation \eqref{kp} becomes
\begin{align}\label{boundary}
 \langle\bar{\omega}, \mathbb{G}(\partial_\eta)\bar{u}\rangle_{L^2(\partial\Omega)} = 0.
\end{align}

Let $\bar{\upsilon}$ be an arbitrary vector on the $\partial \Omega$, It is easy to see that there exists $\bar{u}\in C^{\infty}_0(\Omega)$ such that 
\begin{align}
\mathbb{G}(\partial_\eta)\bar{u} =\bar{\upsilon},
\end{align}
using this in \eqref{boundary}, we have  that

\begin{align}\label{boundary2}
\bar{\omega}|_{\partial\Omega} =\bzero.
\end{align}
From \eqref{cokernel} and \eqref{boundary2} we have the BVP

\begin{align}\label{kk}
P(D)\bar{\omega} = \bzero, \quad \bar{\omega}|_{\partial \Omega} = \bzero.
\end{align}
Equation \eqref{kk} is exactly the equation for solving for the kernel of $P(D)$, therefore, repeating the same argument as before we  conclude that $\bar{\omega}=\bzero.$ Hence co-kernel $\big(P(D)\big)$ = \{$\bzero$\}. Hence the boundary value problem BVP \eqref{BVP4} is uniquely solvable.
\end{proof}


\end{document}